\documentclass[aos,MSNbibl,nameyear,dvips]{arximspdf}
\usepackage{graphicx}

\doi{10.1214/13-AOS1092} 
\volume{41}
\issue{2}
\pubyear{2013}
\firstpage{751}
\lastpage{771}

\makeatletter

\newcommand{\cal}{\mathcal}

\newtheorem{theorem}{Theorem}[section]
\newtheorem{lemma}{Lemma}
\newtheorem{corollary}[theorem]{Corollary}

\newproclaim{Remark}{Remark}

\newcommand{\E}{\mathbb{E}}
\def\E{\mathbb{E}}
\def\cE{{\cal E}}

\renewcommand{\hat}{\widehat}

\makeatother

\begin{document}
\begin{frontmatter}

\title{Density-sensitive semisupervised inference}
\runtitle{Semisupervised inference}

\begin{aug}
\author[A]{\fnms{Martin} \snm{Azizyan}\thanksref{t2}\ead[label=e1]{mazizyan@cs.cmu.edu}},
\author[A]{\fnms{Aarti} \snm{Singh}\thanksref{t2}\ead[label=e2]{aarti@cs.cmu.edu}}
\and
\author[A]{\fnms{Larry} \snm{Wasserman}\corref{}\thanksref{t3}\ead[label=e3]{larry@stat.cmu.edu}}
\runauthor{M. Azizyan, A. Singh and L. Wasserman}
\affiliation{Carnegie Mellon University}
\address[A]{Department of Statistics\\
\quad and Machine Learning Department\\
Carnegie Mellon University\\
Pittsburgh, Pennsylvania 15213\\
USA\\
\printead{e1}\\
\hphantom{E-mail: }\printead*{e2}\\
\hphantom{E-mail: }\printead*{e3}} 
\end{aug}

\thankstext{t2}{Supported by Air Force Grant FA9550-10-1-0382
and NSF Grant IIS-1116458.}

\thankstext{t3}{Supported by NSF Grant DMS-08-06009 and
Air Force Grant FA95500910373.}

\received{\smonth{4} \syear{2012}}
\revised{\smonth{1} \syear{2013}}

%
\begin{abstract}
Semisupervised methods are techniques for using labeled data
$(X_1,Y_1),\ldots,(X_n,Y_n)$ together with unlabeled data
$X_{n+1},\ldots, X_N$ to make predictions. These methods invoke some
assumptions that link the marginal distribution $P_X$ of $X$ to the
regression function $f(x)$. For example, it is common to assume that
$f$ is very smooth over high density regions of $P_X$. Many of the
methods are ad-hoc and have been shown to work in specific examples but
are lacking a theoretical foundation. We provide a minimax framework
for analyzing semisupervised methods. In particular, we study methods
based on metrics that are sensitive to the distribution $P_X$. Our
model includes a parameter $\alpha$ that controls the strength of the
semisupervised assumption. We then use the data to adapt to $\alpha$.
\end{abstract}

%

\end{frontmatter}

\section{Introduction}
\label{sectionintroduction}

Suppose we have data $(X_1,Y_1),\ldots, (X_n,Y_n)$ from a distribution
$P$, where $X_i\in\mathbb{R}^d$ and $Y_i\in\mathbb{R}$. Further, we
have a second set of data $X_{n+1},\ldots, X_N$ from the same
distribution but without the $Y$'s. We refer to ${\cal L} =
\{(X_i,Y_i)\dvtx  i=1,\ldots, n\}$ as the \textit{labeled data} and
${\cal U} = \{X_i\dvtx  i=n+1,\ldots, N\}$ as the \textit{unlabeled data}.
There has been a major effort, mostly in the machine learning
literature, to find ways to use the unlabeled data together with the
labeled data to constuct good predictors of $Y$. These methods are
known as \textit{semisupervised methods}. It is generally assumed that the
$m=N-n$ unobserved labels $Y_{n+1},\ldots, Y_N$ are missing
completely at random and we shall assume this throughout.

To motivate semisupervised inference,
consider the following example.
We download a large number $N$ of webpages $X_i$.
We select a small subset of size $n$ and label these
with some attribute $Y_i$.
The
downloading process is cheap
whereas the labeling process is expensive so typically $N$ is huge
while $n$ is much smaller.

\begin{figure}

\includegraphics{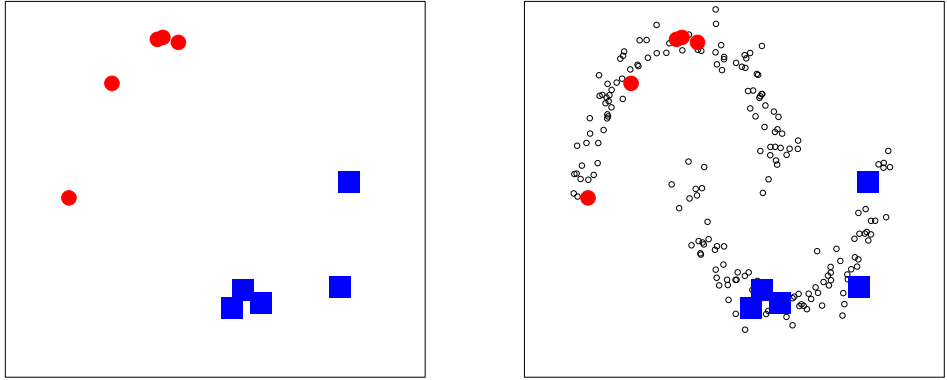}

\caption{The covariate $X=(X_1,X_2)$ is two dimensional.
The response $Y$ is binary and is shown as a square
or a circle. Left: the labeled data.
Right: labeled and unlabeled data.}
\label{figtoy}
\end{figure}

\begin{figure}[b]
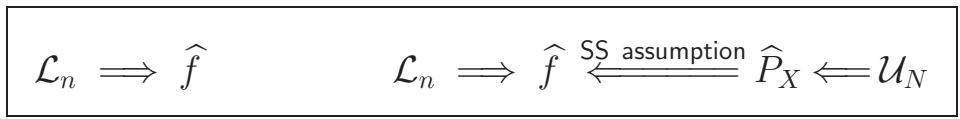

\fbox{\parbox{353pt}{\Large
\[
{\cal L}_n \ \Longrightarrow \ \hat f
\hspace{1in}
{\cal L}_n \ \Longrightarrow \ \hat f\
\stackrel{\mathsf{SS}\ \mathsf{assumption}}{\Longleftarrow \!=\!=\!=\!=}
\hat{P}_X
\Longleftarrow
{\cal U}_N
\]}}
\caption{Supervised learning (left) uses only the labeled data ${\cal L}_n$.
Semisupervised learning (right) uses the unlabeled data ${\cal U}_N$ to estimate
the marginal distribution $P_X$
which helps estimate $f$ if there is some link between $P_X$ and $f$.
This link is the semisupervised (SS) assumption.}
\label{figexplain}
\end{figure}

Figure~\ref{figtoy} shows a toy example
of how unlabeled data can help with prediction.
In this case, $Y$ is binary, $X\in\mathbb{R}^2$
and we want to find the decision boundary
$\{x\dvtx  P(Y=1|X=x)=1/2\}$.
The left plot shows a few labeled data points
from which it would be challenging to find the boundary.
The right plot shows labeled and unlabeled points.
The unlabeled data show that there are two clusters.
If we make the seemingly reasonable assumption
that $f(x)=P(Y=1|X=x)$ is very smooth over the two clusters,
then identifying the decision boundary becomes much easier.
In other words, if we assume some link between $P_X$ and $f$, then we
can use the unlabeled data;
see Figure~\ref{figexplain}.

The assumption that
the regression function $f(x) = \mathbb{E}(Y|X=x)$ is very smooth over
the clusters
is known as the \textit{cluster assumption}.
In the special case where the clusters are low-dimensional submanifolds,
the assumption is called
the \textit{manifold assumption}.
These assumptions link the regression function $f$ to the distribution
$P_X$ of $X$.

Many semisupervised methods are developed based on the above
assumptions, although this is not always made explicit. Even with such
a link, it is not obvious that semisupervised methods
will outperform supervised methods.
Making precise
how and when these assumptions actually improve inferences is
surprisingly elusive, and most papers
do not address this issue; some exceptions are \citet{rigollet07},
\citet{SSLTR}, \citet{LWnips07}, \citet{nadler09},
\citet{ben-davidcolt08}, \citet{NIPS20091025},
\citet{belkinniyogi} and \citet{partha}.
These authors have shown that the degree to which
unlabeled data improves performance is very sensitive to the cluster
and manifold assumptions. In this paper, we introduce \textit{adaptive
semisupervised inference}. We define a parameter $\alpha$ that
controls the sensitivity of the distance metric to the density, and
hence the strength of the semisupervised assumption. When $\alpha= 0$
there is no semisupervised assumption, that is, there is no link
between $f$ and $P_X$. When $\alpha= \infty$ there is a very strong
semisupervised assumption. We use the data to estimate $\alpha$, and
hence we adapt to the appropriate assumption linking $f$ and $P_X$.
In addition, we should add that we focus on regression while most
previous literature
only deals with binary outcomes (classification).

This paper makes the following contributions:
\begin{longlist}[(6)]
\item[(1)] We formalize the link between
the regression function $f$ and the marginal distribution
$P_X$ by defining a class of function spaces
based on a metric that depends on $P_X$.
This is called a \textit{density sensitive metric}.
\item[(2)] We show how to consistently estimate the density-sensitive metric.
\item[(3)] We propose a semi-supervised kernel estimator based on
the density-sensitive metric.
\item[(4)]
We provide some minimax bounds and
show that under some conditions the semisupervised method
has smaller predictive risk than any supervised method.
\item[(5)] The function classes depend on a parameter $\alpha$ that
controls how strong the
semisupervised assumption is.
We show that it is possible to adapt to $\alpha$.
\item[(6)] We provide numerical simulations to support the theory.
\end{longlist}

We now give an informal statement of our main results.
In Section~\ref{sectionminimax} we define a nonparametric class of
distributions
${\cal P}_n$.
Let $0 < \xi< d-3$ and assume that
$m \geq n^{2/(2+\xi)}$.
Let ${\cal S}_n$ denote the set of supervised estimators;
these estimators use only the labeled data.
Let ${\cal SS}_N$ denote the set of semisupervised estimators;
these estimators use the labeled data and unlabeled data.
Then:
\begin{longlist}[(3)]
\item[(1)]
(Theorem~\ref{thmupper-bound} and Corollary~\ref{corollaryupper}.)
There is a semisupervised estimator $\hat f$ such that
%
\begin{equation}
\sup_{P\in{\cal P}_n}R_P(\hat f) \leq \biggl(
\frac{C}{n} \biggr)^{{2}/({2+\xi})},
\end{equation}
where $R_P(\hat f)$ is the risk of the estimator $\hat f$ under
distribution $P$.
\item[(2)] (Theorem~\ref{thmlower-bound}.)
For supervised estimators ${\cal S}_n$ we have
%
\begin{equation}
\inf_{\hat f\in{\cal S}_n}\sup_{P\in{\cal P}_n} R_P(\hat
f) \geq \biggl(\frac{C}{n} \biggr)^{{2}/({d-1})}.
\end{equation}
\item[(3)] Combining these two results we conclude that
%
\begin{equation}
\frac{\inf_{\hat f\in{\cal SS}_N}\sup_{P\in{\cal P}_n} R_P(\hat f)} {
\inf_{\hat f\in{\cal S}_n}\sup_{P\in{\cal P}_n} R_P(\hat f)} \leq
\biggl(\frac{C}{n} \biggr)^{{2(d-3-\xi)}/({(2+\xi)(d-1)})} \to0
\end{equation}
and hence, semisupervised estimation dominates supervised estimation.
\end{longlist}

\begin{Remark*}
We assume, as is standard in the literature
on semisupervised learning,
that the margial $P_X$ is the same for the labeled and unlabeled
data. Extensions to the case where the marginal distribution changes
are possible, but are beyond the scope of the paper.
\end{Remark*}

\textit{Related work.}
There are a number of papers that
discuss conditions under which
semisupervised methods can succeed
or that discuss metrics that are useful for
semisupervised methods.
These include Castelli and Cover
(\citeyear{castellicover95,castellicover96}),
\citet{ratsaby1995learning},
\citet{bousquet04}, \citet{SSLTR}, \citet{LWnips07},
\citet{NIPS20091025},
\citet{ben-davidcolt08},
\citet{nadler09}, \citet{orlitsky}, \citet{bigral11},
\citet{belkinniyogi}, \citet{partha}
and references therein.
Papers on semisupervised inference in the statistics
literature are rare; some exceptions include
\citet{Culp2008}, \citet{Culp2011} and
\citet{West2007}.
To the best of our knowledge,
there are no papers that explicitly study
adaptive methods that allow the data to choose the
strength of the semisupervised assumption.

There is a connection between our work on the semisupervised
classification method in \citet{rigollet07}. He divides the covariate
space ${\cal X}$ into clusters $C_1,\ldots, C_k$ defined by the upper
level sets $\{p_X > \lambda\}$ of the density $p_X$ of $P_X$. He
assumes that the indicator function $I(x) = I(p(y|x) > 1/2)$ is
constant over each cluster $C_j$. In our regression framework, we could
similarly assume that
\[
f(x) = \sum_{j=1}^k f_{\theta_j}(x)
I(x\in C_j) + g(x) I(x\in C_0),
\]
where $f_\theta(x)$ is a parametric regression function, $g$ is a
smooth (but nonparametric function) and $C_0 = {\cal X} -
\bigcup_{j=1}^k C_j$. This yields\vspace*{1pt} parametric,
dimension-free rates over ${\cal X}-C_0$. However, this creates a
rather unnatural and harsh boundary at $\{x\dvtx  p_X(x) = \lambda\}$.
Also, this does not yield improved rates over $C_0$. Our approach may
be seen as a smoother version of this idea.\vadjust{\goodbreak}

\textit{Outline.} This paper is organized as follows. In Section
\ref{secsetup} we give definitions and assumptions. In Section
\ref{sectiondefine-metrics} we define density sensitive metrics and the
function spaces defined by these metrics. In Section~\ref{secest} we
define a density sensitive semisupervised estimator, and we bound its
risk. In Section~\ref{sectionminimax} we present some minimax results.
We discuss adaptation in Section~\ref{secadap}. We provide simulations
in Section~\ref{secsim}. Section~\ref{secdisc} contains the closing
discussion. Many technical details and extensions are contained in the
supplemental article
[\citet{AzizyanSSLsupplement}].

\section{Definitions}
\label{secsetup}

Recall that $X_i\in\mathbb{R}^d$ and $Y_i\in\mathbb{R}$.
Let
%
\begin{equation}
{\cal L}_n = \bigl\{(X_1,Y_1),\ldots,(X_n,Y_n)\bigr\}
\end{equation}
be an i.i.d. sample from $P$. Let $P_X$
denote the $X$-marginal of $P$, and let
%
\begin{equation}
{\cal U}_N=\{X_{n+1},\ldots, X_N\}
\end{equation}
be an i.i.d. sample from $P_X$.

Let $f(x) \equiv f_P(x)= \mathbb{E}(Y|X=x)$. An estimator of $f$ that
is a function of ${\cal L}_n$ is called a \textit{supervised learner},
and the set of such estimators is denoted by ${\cal S}_n$. An estimator
that is a function of ${\cal L}_n \cup{\cal U}_N$ is called a
\textit{semisupervised learner}, and the set of such estimators is
denoted by ${\cal SS}_N$. Define the risk of an estimator $\hat f$ by
%
\begin{equation}
R_P(\hat f) = \mathbb{E}_P \biggl[\int\bigl(\hat f(x)
- f_P(x)\bigr)^2 \,dP(x) \biggr],
\end{equation}
where $\mathbb{E}_P$ denotes the expectation over data drawn
from the distribution $P$.
Of course,
${\cal S}_n \subset{\cal SS}_N$ and
hence
\[
\inf_{\hat g\in{\cal SS}_N}\sup_{P\in{\cal P}}R_P(\hat g)
\leq \inf_{\hat g\in{\cal S}_n}\sup_{P\in{\cal P}}R_P(
\hat g).
\]
We will show that,
under certain conditions,
semisupervised methods outperform
supervised methods in the sense that the left-hand side of the above equation
is substantially smaller than the right-hand side.
More precisely,
for certain classes of
distributions
${\cal P}_n$,
we show that
%
\begin{equation}
\frac{\inf_{\hat g\in{\cal SS}_N}\sup_{P\in{\cal P}_n}R_P(\hat g)} {
\inf_{\hat g\in{\cal S}_n}\sup_{P\in{\cal P}_n}R_P(\hat g)} \to0
\end{equation}
as $n\to\infty$.
In this case we say that
semisupervised learning is \textit{effective}.

\begin{Remark*}
In order for the asymptotic analysis to reflect the
behavior of finite samples, we need to let ${\cal P}_n$ to change with
$n$, and we need $N = N(n)\to\infty$ and $n/N(n)\to0$ as $n\to\infty$.
As an analogy, one needs to let the number of covariates in a
regression problem increase with the sample size to develop relevant
asymptotics for high-dimensional regression. Moreover, ${\cal P}_n$
must have distributions that get more concentrated as $n$ increases.
The reason is that if $n$ is very large and $P_X$ is smooth, then there
is no advantage to semisupervised inference. This is consistent with
the finding in \citet{ben-davidcolt08} who show that if $P_X$ is smooth,
then ``$\ldots$ knowledge of that distribution cannot improve the
labeled sample complexity by more than a constant
factor.''
\end{Remark*}

\textit{Other notation.}
If $A$ is a set and $\delta\geq0$, we define
\[
A\oplus\delta= \bigcup_{x\in A} B(x,\delta),
\]
where $B(x,\delta)$ denotes a ball of radius $\delta$ centered at $x$.
Given a set $A\subseteq\mathbb{R}^d$, define
$d_A(x_1,x_2)$
to be the length of the shortest path in $A$ connecting $x_1$ and $x_2$.

We write
$a_n = O(b_n)$ if
$|a_n/b_n|$ is bounded for all large $n$.
Similarly,
$a_n = \Omega(b_n)$ if
$|a_n/b_n|$ is bounded away from 0 for all large $n$.
We write
$a_n \asymp b_n$ if
$a_n = O(b_n)$ and
$a_n = \Omega(b_n)$.
We also write
$a_n \preceq b_n$ if
there exists $C>0$ such that
$a_n \leq C b_n$ for all large $n$.
Define
$a_n \succeq b_n$ similarly.
We use symbols of the form
$c,c_1,c_2,\ldots, C,C_1,C_2,\ldots$
to denote generic positive constants whose
value can change in different expressions.

\section{Density-sensitive function spaces}
\label{sectiondefine-metrics}

We define a smoothed version of $P_X$ as follows.
(This is needed since we allow the marginal distribution $P_X$ to be singular.)
Let $K$ denote a symmetric kernel on $\mathbb{R}^d$
with compact support,
let $\sigma>0$
and define
%
\begin{equation}
p_\sigma(x) \equiv p_{X,\sigma}(x) = \int\frac{1}{\sigma^d}K \biggl(
\frac{\|x-u\|}{\sigma} \biggr) \,dP_X(x).
\end{equation}
Thus,
$p_{X,\sigma}$ is the density of the convolution
$P_{X,\sigma}=P_X\star\mathbb{K}_\sigma$
where
$\mathbb{K}_\sigma$ is the measure with density
$K_\sigma(\cdot) = \sigma^{-d}K(\cdot/\sigma)$.
$P_{X,\sigma}$ always has a density even if $P_X$ does not.
This is important because, in high-dimensional problems,
it is not uncommon to find that $P_X$ can be highly concentrated near a
low-dimensional manifold. These are
precisely the cases where semisupervised methods are often useful
[\citet{ben-davidcolt08}].
Indeed, this was one of the original motivations for semisupervised inference.
We define $P_{X,0} = P_X$.
For notational simplicity, we shall sometimes drop the $X$ and simply write
$p_\sigma$ instead of $p_{X,\sigma}$.

\subsection{The exponential metric}

\begin{figure}

\includegraphics{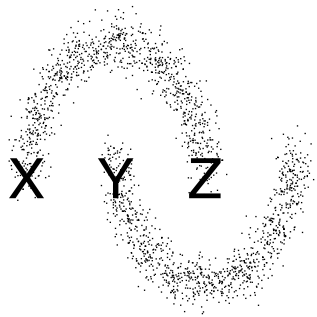}

\caption{With a density metric, the points $X$ and $Z$ are
closer than the points $X$ and $Y$ because there is a high density
path connecting $X$ and $Z$.}
\label{figdensitymetricpicture}
\end{figure}

Following previous work in the area, we will assume that the regression
function is smooth in regions where $P_X$ puts lots of mass. To make
this precise, we define a \textit{density sensitive metric} as follows.
For any pair $x_1$ and $x_2$ let $\Gamma(x_1,x_2)$ denote the set of
all continuous finite curves from $x_1$ to $x_2$ with unit speed
everywhere, and let $L(\gamma)$ be the length\vadjust{\goodbreak} of curve $\gamma$; hence
$\gamma(L(\gamma))=x_2$. For any $\alpha\geq0$ define the
\textit{exponential metric}
%
\begin{equation}
D(x_1,x_2) \equiv D_{P,\alpha,\sigma}(x_1,x_2)
= \inf_{\gamma\in\Gamma(x_1,x_2)} \int_0^{L(\gamma)} \exp
\bigl[-\alpha p_{X,\sigma
}\bigl(\gamma(t)\bigr) \bigr]\,dt.
\end{equation}
In the supplement, we also
consider a second metric, the \textit{reciprocal metric}.
Large $\alpha$ makes points connected by high density
paths closer; see Figure~\ref{figdensitymetricpicture}.
%
Note that $\alpha=0$ corresponds to Euclidean distance.
Similar definitions are used in \citet{orlitsky}, \citet{bigral11} and
\citet{bousquet04}.

\subsection{The regression function}

Recall that $f(x)\equiv f_P(x)=E(Y|X=x)$ denotes the regression
function. We assume that $X\in[0,1]^d\equiv{\cal X}$ and that $|Y| \leq
M$ for some finite constant $M$.\setcounter{footnote}{2}\footnote{The
results can be extended to unbounded $Y$ with suitable conditions on
the tails of the distribution of $Y$.} We formalize the semisupervised
smoothness assumption by defining the following scale of function
spaces. Let ${\cal F} \equiv{\cal F}(P,\alpha,\sigma,L)$ denote the set
functions $f\dvtx  [0,1]^d \to\mathbb{R}$ such that, for all
$x_1,x_2\in{\cal X}$,
%
\begin{equation}
\bigl|f(x_1) - f(x_2)\bigr| \leq L D_{P,\alpha,\sigma}(x_1,x_2).
\end{equation}
Let
${\cal P}(\alpha,\sigma,L)$ denote all joint distributions for $(X,Y)$
such that $f_P\in{\cal F}(P,\alpha,\break\sigma,L)$
and such that $P_X$ is supported on ${\cal X}$.

\subsection{Properties of the function spaces}

Let $B_{P,\alpha,\sigma}(x,\varepsilon)=\{z\dvtx  D_{P,\alpha,\sigma
}(x,\break z) \leq\varepsilon\}$
be a ball of size $\varepsilon$.
Let $S_P$ denote the support of $P$, and
let ${\cal N}_{P,\alpha,\sigma}(\varepsilon)$ denote the covering number,
the smallest number of balls of size $\varepsilon$
required to cover $S_P$.
The covering number measures the size of the function space,
and the variance of any regression estimator
on the space ${\cal F}(P,\alpha,\sigma,L)$ depends on this covering number.
Here, we mention a few properties of
${\cal N}_{P,\alpha,\sigma}(\varepsilon)$.

In the Euclidean case $\alpha=0$,
we have
${\cal N}_{P,0,\sigma}(\varepsilon) \leq(C/\varepsilon)^d$.
But when $\alpha>0$ and $P$ is concentrated on or near
a set of dimension less than $d$, the
${\cal N}_{P,\alpha,\sigma}(\varepsilon)$ can be much smaller than
$(C/\varepsilon)^d$.
The next result gives a few examples
showing that concentrated distributions have small covering numbers.
We say that a set $A$ is \textit{regular} if there is a $C>0$ such that,
for all small $\varepsilon>0$,
%
\begin{equation}
\sup_{ \stackrel{x,y\in A}{\|x-y\| \leq\varepsilon}} \frac
{d_A(x,y)}{\|x-y\|} \leq C,
\end{equation}
where $d_A(x_1,x_2)$
is the length of the shortest path in $A$ connecting $x_1$ and $x_2$.
Recall that $S_P$ denotes the support of $P$.

\begin{lemma}
Suppose that $S_P$ is regular.
\begin{longlist}[(4)]
\item[(1)] For all $\alpha$, $\sigma$ and $P$,
${\cal N}_{P,\alpha,\sigma}(\varepsilon)\preceq\varepsilon^{-d}$.
\item[(2)] Suppose that $P = \sum_{j=1}^k \delta_{x_j}$
where $\delta_x$ is a point mass at $x$.
Then, for any $\alpha\geq0$ and any $\varepsilon>0$,
${\cal N}_{P,\alpha,\sigma}(\varepsilon)\leq k$.
\item[(3)] Suppose that $\mathsf{dim}(S_P) = r < d$.
Then,
${\cal N}_{P,\alpha,\sigma}(\varepsilon) \preceq\varepsilon^{-r}$.
\item[(4)]
Suppose that $S_P = W\oplus\gamma$ where
$\mathsf{dim}(W) =r<d$.
Then, for $\varepsilon\geq C \gamma$,
${\cal N}_{P,\alpha,\sigma}(\varepsilon) \preceq
(\frac{1}{\varepsilon} )^{r}$.
\end{longlist}
\end{lemma}

\begin{pf}
(1) The first statement follows since the covering number of $S_P$ is
no more than the covering number of $[0,1]^d$ and on $[0,1]^d$,
$D_{P,\alpha,\sigma}(x,y) \leq\|x-y\|$. Now $[0,1]^d$ can be covered
$O(\varepsilon^{-d})$ Euclidean balls.

(2) The second statement follows since $\{\{x_1\},\ldots, \{x_k\}\}$
forms an $\varepsilon$-covering for any $\varepsilon$.

(3) We have that $D_{P,\alpha,\sigma}(x,y) \leq d_{S_P}(x,y)$.
Regularity implies that, for small $d_{S_P}(x,y)$,
$D_{P,\alpha,\sigma}(x,y) \leq c\|x-y\|$. We can thus cover $S_P$ by
$C\varepsilon^{-r}$ balls of size~$\varepsilon$.

(4) As in (3), cover $W$ with $N=O(\varepsilon^{-r})$ balls of $D$ size
$\varepsilon$. Denote these balls by $B_1,\ldots, B_N$. Define $C_j =
\{x\in S_P\dvtx  d_{S_P}(x,B_j) \leq\gamma\}$. The $C_j$ form a
covering of size $N$ and each $C_j$ has $D_{P,\alpha,\sigma}$ diameter
$\max\{\varepsilon,\gamma\}$.
\end{pf}

\section{Semisupervised kernel estimator}
\label{secest}

We consider the following
semisupervised
estimator which uses a kernel that is
sensitive to the density.
Let $Q$ be a kernel and let
$Q_h(x) = h^{-d} Q(x/h)$.
Let
%
\begin{equation}
\label{eqdefine-estimator} \hat f_{h,\alpha,\sigma}(x) = \frac{\sum^n_{i=1}Y_i   Q_h
(\widehat D_{\alpha,\sigma}(x,X_i) )} {
\sum^n_{i=1}Q_h (\widehat D_{\alpha,\sigma}(x,X_i) )},
\end{equation}
where
%
\begin{eqnarray}
\hat D_{\alpha,\sigma}(x_1,x_2) &=& \inf
_{\gamma\in\Gamma(x_1,x_2)} \int_0^{L(\gamma)}\exp \bigl[ -
\alpha\hat p_\sigma\bigl(\gamma (t)\bigr) \bigr] \,dt,
\\
%
\hat{p}_\sigma(x) &=& \frac{1}{m}\sum_{i=1}^m
\frac{1}{\sigma^d} K \biggl( \frac{\|x-X_{i+n}\|}{\sigma} \biggr),
\end{eqnarray}
and
$m = N-n$ denotes the number of unlabeled points.
We use a kernel estimator for the regression function
because it is simple, commonly used and, as we shall see, has
a fast rate of convergence in the semisupervised case.

The estimator
$\hat D_{\alpha,\sigma}(x_1,x_2)$
is discussed in detail in the supplement
where we study its properties
and we give an algorithm for
computing it.

Now we give an upper bound
on the risk of $\hat f_{h,\alpha,\sigma}$.
In the following we take, for simplicity,
$Q(x) = I(\|x\| \leq1)$.

\begin{theorem}
\label{thmupper-bound}
Suppose that $|Y|\leq M$.
Define the event $\mathcal{G}_m = \{ \|\hat p_\sigma- p_\sigma
\|_\infty\leq\varepsilon_m\}$
(which depends on the unlabeled data)
and suppose that
$\mathbb{P}({\cal G}_m^c)\leq1/m$.
Then,
for every $P\in{\cal P}(\alpha,\sigma,L)$,
%
\begin{equation}
R_P(\hat f_{h,\alpha,\sigma})\leq L^2 \bigl(h
e^{\alpha\varepsilon_m}\bigr)^{2} + \frac{M^2 (2 + {1}/{e} )  {\cal
N}(P,\alpha,\sigma,e^{-\varepsilon_m\alpha}h/2)}{n} +
\frac{4M^2}{m}.\hspace*{-32pt}
\end{equation}
\end{theorem}

\begin{pf}
The risk is
\begin{eqnarray*}
R_P(\hat f) &=& \mathbb{E}_{n,N} \biggl[ (1- \mathcal{G}_m)
\int\bigl(\hat{f}_{h,\alpha,\sigma}(x) - f(x) \bigr)^2 \,dP(x) \biggr]\\
&&{} +
\mathbb{E}_{n,N} \biggl[ \mathcal{G}_m
\int\bigl(\hat{f}_{h,\alpha,\sigma}(x) - f(x) \bigr)^2 \,dP(x) \biggr].
\end{eqnarray*}
Since $|Y| \leq M$
and $\sup_x|\hat f(x)| \leq M$,
\[
\mathbb{E}_{n,N} \biggl[ (1-\mathcal{G}_m) \int\bigl(
\hat{f}_{h,\alpha,\sigma}(x) - f(x)\bigr)^2 \,dP(x) \biggr]\leq 4
M^2 \mathbb{P}\bigl(\mathcal{G}_m^c\bigr)
\leq\frac{4M^2}{m}.
\]
Now we bound the second term.

Condition on the unlabeled data.
Replacing the Euclidean distance with $\hat D_{\alpha,\sigma}$ in
the proof of Theorem 5.2 in \citet{gyorfi2002nonparametric},
we have that
\begin{eqnarray*}
&&
\mathbb{E}_{n} \biggl[ \int\bigl(\hat{f}_{h,\alpha,\sigma}(x) - f(x)
\bigr)^2 \,dP(x) \biggr]\\
&&\qquad\leq L^2 R^{2} +
\frac{M^2 (2 + {1}/{e} ) \int{dP(x)}/{P(\hat
B_{\alpha,\sigma}(x,h))}}{n},
\end{eqnarray*}
where
\[
R = \sup \bigl\{ D_{P,\alpha,\sigma}(x_1,x_2)\dvtx
(x_1,x_2) \mbox{ such that }\hat D_{\alpha,\sigma}(x_1,x_2)
\leq h \bigr\}
\]
and
$\hat B_{\alpha,\sigma}(x,h) = \{z\dvtx  \hat{D}_{\alpha,\sigma
}(x,z)\leq h\}$.
On the event ${\cal G}_m$,
we have from Lemma~2 in the supplement 
that
$e^{-\alpha\varepsilon_m} D_{\alpha,\sigma}(x_1,x_2) \leq
\hat D_{\alpha,\sigma}(x_1,x_2) \leq
e^{\alpha\varepsilon_m} D_{\alpha,\sigma}(x_1,x_2)$
for all $x_1,x_2$.
Hence,
$R^{2}\leq e^{2\alpha\varepsilon_m} h^{2}$ and
\[
\int\frac{dP(x)}{P(\hat B_{\alpha,\sigma}(x,h))}\leq \int\frac{dP(x)}{P(B_{P,\alpha,\sigma}(x, e^{-\alpha\varepsilon_m}h))}.
\]
A simple covering argument [see page 76 of \citet{gyorfi2002nonparametric}] shows that,
for any $\delta>0$,
\[
\int\frac{dP(x)}{P(B_{P,\alpha,\sigma}(x,\delta))} \leq{\cal N}(P,\alpha,\sigma,\delta/2).
\]
The result follows.
\end{pf}

\begin{corollary}
\label{corollaryupper}
If
${\cal N}(P,\alpha,\sigma,\delta) \leq(C/\delta)^\xi$
for
$\delta\geq(1/2)e^{-\alpha\varepsilon_m} (n\times \break e^{2\alpha\varepsilon
_m})^{-{1}/({2+\xi})}$
and $N \geq2n$,
then
%
\begin{equation}
R_P(\hat f_{\alpha,\sigma,h}) \leq e^{\alpha\varepsilon_m (2\vee\xi)} \biggl[
L^2 h^2 + \frac{1}{n} \biggl(\frac{C}{h}
\biggr)^\xi \biggr] + \frac{4 M^2}{m}.
\end{equation}
Hence, if
$m \geq n^{2/(2+\xi)}$
and $h \asymp(n e^{\alpha\varepsilon_m (2-\xi)})^{-{1}/({2+\xi})}$,
then
%
\begin{equation}
\sup_{P\in{\cal P}(\alpha,\sigma,L)}R_P(\hat f_{h,\alpha,\sigma}) \preceq
\biggl(\frac{C}{n} \biggr)^{{2}/({2+\xi})}.
\end{equation}
\end{corollary}

\section{Minimax bounds}
\label{sectionminimax}

To characterize when semisupervised methods
outperform supervised methods,
we show that there is a class
of distributions~${\cal P}_n$
(which we allow to change with $n$) such that
$R_{SS}$ is much smaller than
$R_{S}$, where
\[
R_S = \inf_{\hat f\in{\cal S}_n}\sup_{P\in{\cal P}_n}
R_P(\hat f) \quad\mbox{and}\quad R_{SS}= \inf
_{\hat f\in{\cal SS}_N}\sup_{P\in{\cal P}_n} R_P(\hat f).
\]
To do so, it suffices to find a lower bound on
$R_{S}$ and an upper bound on $R_{SS}$.
Intuitively, ${\cal P}_n$ should be a set
distributions whose $X$-marginals are highly concentrated
on or near lower-dimensional sets, since this is where semisuspervised
methods deliver improved performance.
Indeed, as we mentioned earlier,
for very smooth distributions $P_X$
we do not expect semisupervised learners to offer much improvement.

\subsection{The class ${\cal P}_n$}

Here we define the class ${\cal P}_n$.
Let $N= N(n)$ and $m = m(n)=N-n$ and define
%
\begin{equation}
\varepsilon_m \equiv\varepsilon(m,\sigma) = \sqrt{\frac{C \log m}{m
\sigma^d}}.
\end{equation}
Let $\xi\in[0,d-3)$, $\gamma>0$ and define
%
\begin{equation}
{\cal P}_{n} = \bigcup_{(\alpha,\sigma)\in{\cal A}_n\times\Sigma_n} {\cal Q}(
\alpha,\sigma,L),
\end{equation}
where
${\cal Q}(\alpha,\sigma,L)\subset{\cal P}(\alpha,\sigma,L)$ and
${\cal A}_n\times\Sigma_n \subset[0,\infty]^2$
satisfy the following conditions:
\begin{eqnarray*}
&&\mbox{(C1)}\quad {\cal Q}(\alpha,\sigma, L) \\
&&\hphantom{\mbox{(C1)}}\quad\qquad= \biggl\{P\in{
\cal P}(\alpha,\sigma,L)\dvtx  {\cal N}(P,\alpha,\sigma,\varepsilon) \leq \biggl(
\frac{C}{\varepsilon} \biggr)^\xi\ \forall \varepsilon\geq \biggl(
\frac{1}{n} \biggr)^{{1}/({2+\xi})} \biggr\};
\\
&&\mbox{(C2)}\quad \alpha\leq\frac{\log2}{\varepsilon(m,\sigma)};
\\
&&\mbox{(C3)}\quad \biggl(\frac{1}{m} \biggr)^{{1}/({d(1+\gamma)})} \leq \sigma\leq
\frac{1}{4C_0} \biggl(\frac{1}{n} \biggr)^{{1}/({d-1})},
\end{eqnarray*}
where $C_0$ is the diameter of the support of $K$.

Here are some remarks about ${\cal P}_{n}$:
\begin{longlist}[(3)]
\item[(1)] (C2) implies that
$e^{\alpha\varepsilon_m} \leq2$ and hence,
(C3) and Theorem 1.3 in the supplement 
$(1/2)D_{P,\alpha,\sigma}(x_1,x_2) \leq
\hat D_{\alpha,\sigma}(x_1,x_2) \leq
2D_{P,\alpha,\sigma}(x_1,x_2)$
with probability at least $1-1/m$.
%
\item[(2)] The constraint in (C1) on ${\cal N}(\varepsilon)$ holds
whenever $P$ is concentrated on or near a
set of dimension less than $d$ and $\alpha/\sigma^d$ is large.
The constraint
does not need to hold for
arbitrarily small $\varepsilon$.
\item[(3)] Some papers on semisupervised learning
simply assume that $N=\infty$
since in practice $N$ is usually very large compared to $n$.
In that case,
there is no upper bound on $\alpha$ and no
lower bound on $\sigma$.
\end{longlist}

The class ${\cal P}_n$ may seem complicated.
This is because showing conditions where
semisupervised learning provably outperforms
supervised learning is subtle.
Intuitively, the class ${\cal P}_n$
is simply the set of high concentrated distributions
with $\alpha/\sigma$ large.

\subsection{Supervised lower bound}

\begin{theorem}
\label{thmlower-bound}
Suppose that
$m \geq n^{{d(1+\gamma)}/({d-1})}$.
There exists $C>0$ such that
%
\begin{equation}
R_S = \inf_{\hat f\in{\cal S}_n}\sup_{P\in{\cal P}_{n}}
R_P(\hat f) \geq \biggl(\frac{C}{n} \biggr)^{{2}/({d-1})}.
\end{equation}
\end{theorem}

\begin{pf}
Let
$A_1$ and $A_0$ be the top and bottom of the cube ${\cal X}$,
\begin{eqnarray*}
A_1 &=& \bigl\{ (x_1,\ldots, x_{d-1},1)\dvtx  0
\leq x_1,\ldots, x_{d-1} \leq 1\bigr\},
\\
A_0 &=& \bigl\{ (x_1,\ldots, x_{d-1},0)\dvtx  0
\leq x_1,\ldots, x_{d-1} \leq 1\bigr\}.
\end{eqnarray*}
Fix $\varepsilon= n^{-{1}/({d-1})}$.
Let $q = (1/\varepsilon)^{d-1}\asymp n$.
For any integers
$s=(s_1,\ldots, s_{d-1})\in N^{d-1}$
with $0 \leq s_i \leq1/\varepsilon$,
define the tendril
\[
\bigl\{ (s_1 \varepsilon,s_2 \varepsilon,\ldots,
s_{d-1}\varepsilon,x_d)\dvtx  \varepsilon\leq x_d \leq1-
\varepsilon\bigr\}.\vadjust{\goodbreak}
\]
There are
$q = (1/\varepsilon)^{d-1}\approx n$
such tendrils.
Let us label the tendrils as
$T_1,\ldots, T_q$.
Note that the tendrils do not quite join up with $A_0$ or $A_1$.

Let
\[
C = A_0 \cup A_1 \cup \Biggl(\bigcup
_{j=1}^q T_j \Biggr).
\]
Define a measure $\mu$ on $C$ as follows:
\[
\mu= \frac{1}{4} \mu_0 + \frac{1}{4}
\mu_1 + \frac{1}{2 q (1-2\varepsilon)} \sum_j
\nu_j,
\]
where $\mu_0$ is $(d-1)$-dimensional Lebesgue measure on $A_0$,
$\mu_1$ is $(d-1)$-dimensional Lebesgue measure on $A_1$ and
$\nu_j$ is one-dimensional Lebesgue measure on $T_j$.
Thus, $\mu$ is a probability measure and
$\mu(C)=1$.

\begin{figure}

\includegraphics{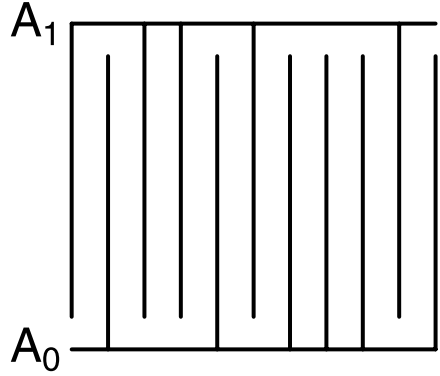}

\caption{The extended tendrils used in the proof of the lower bound,
in the special case where $d=2$.
Each tendril has length $1-\varepsilon$ and joins up with either the top $A_1$
or bottom $A_0$ but not both.}
\label{figsimpletendrils}
\end{figure}

Now we define extended tendrils that are joined to the top or bottom of
the cube
(but not both).
See Figure~\ref{figsimpletendrils}.
If
\[
T_j = \bigl\{ (s_1 \varepsilon,s_2 \varepsilon,\ldots, s_{d-1}\varepsilon,x_d)\dvtx  \varepsilon\leq x_d
\leq1-\varepsilon\bigr\}
\]
is a tendril, define
its extensions
\begin{eqnarray*}
T_{j,0} &=& \bigl\{ (s_1 \varepsilon,s_2
\varepsilon,\ldots, s_{d-1}\varepsilon,x_d)\dvtx  0 \leq
x_d \leq1-\varepsilon\bigr\},
\\
T_{1,j} &=& \bigl\{ (s_1 \varepsilon,s_2
\varepsilon,\ldots, s_{d-1}\varepsilon,x_d)\dvtx  \varepsilon\leq
x_d \leq1\bigr\}.
\end{eqnarray*}
Given
$\omega\in\Omega=\{0,1\}^q$, let
\[
S_\omega= A_0 \cup A_1 \cup \Biggl(\bigcup
_{j=1}^q T_{j,\omega
_j} \Biggr)
\]
and
\[
P_{\omega,X} = \frac{1}{4}\mu_0 + \frac{1}{4}
\mu_1 + \frac{1}{2 q(1-\varepsilon)}\sum_j
\nu_{j,\omega_j},
\]
where
$\nu_{j,\omega_j}$ is one-dimensional Lebesgue measure on
$T_{j,\omega_j}$.
This
$P_{\omega,X}$ is a probability measure supported on $S_\omega$.

Notice that $S_\omega$ consists of two connected components, namely,
\[
U_{\omega}^{(1)} = A_1 \cup \biggl(\bigcup
_{j: \omega_j=1} T_{j,\omega_j} \biggr) \quad\mbox{and}\quad
U_{\omega}^{(0)} = A_0 \cup \biggl(\bigcup
_{j: \omega_j=0} T_{j,\omega_j} \biggr).
\]
Let
\[
f_\omega(x) = \frac{L \varepsilon}{8} I\bigl(x\in U_{\omega}^{(1)}
\bigr).
\]
Finally,
we define
$P_\omega= P_{\omega,X}\times P_{\omega, Y|X}$
where
$P_{\omega, Y|X}$
is a point mass at $f_\omega(X)$.
Define
$d^2(f,g) = \int(f(x)-g(x))^2 \,d\mu(x)$.

We complete the proof with a series of claims.\vspace*{9pt}

\textit{Claim} 1: For each $\omega\in\Omega$, $P_\omega\in{\cal
P}_{n}$.\vspace*{9pt}

\textit{Proof}:
Let
\[
\sigma= \biggl(\frac{1}{m} \biggr)^{{1}/({d(1+\gamma)})}
\]
and let
%
\begin{equation}
\label{eqalp} \frac{3}{2+\xi} \frac{\log m}{m^{{1}/({1+\gamma})}} \leq\alpha \leq \sqrt{
\frac{m^{{\gamma}/({1+\gamma})}}{\log m}}.
\end{equation}
It follows that (C2) and (C3) hold. We must verify (C1). If $x$ and $y$
are in the same connected component, then
$|f_\omega(x)-f_\omega(y)|=0$. Now let $x$ and $y$ be in different
components, that is, $x\in U_\omega^{(1)}, y\in U_\omega^{(0)}$. Let us
choose $x$ and $y$ as close as possible in Euclidean distance; hence
$\|x-y\|=\varepsilon$. Let $\gamma$ be any path connecting $x$ to $y$.
Since $x$ and $y$ lie on different components, there exists a subset
$\gamma_0$ of $\gamma$ of length at least $\varepsilon$ on which
$P_\omega$ puts zero mass. By assumption (C3),
$\sigma\leq\varepsilon/(4C_0)$ and hence $P_{X,\sigma}$ puts zero mass
on the portion of $\gamma_0$ that is at least $C_0\sigma$ away from the
support of $P_\omega$. This has length at least $\varepsilon- 2 C_0
\sigma\geq\varepsilon/2$. Since $p_{X,\sigma}(x)=0$ on a portion
of~$\gamma_0$,
\[
D_{P,\alpha,\sigma}(x,y) \geq\frac{\varepsilon}{2} = \frac{\|x-y\|}{2}.
\]
Hence,
$\|x-y\| \leq2 D_{P,\alpha,\sigma}(x,y)$.
Then
\[
\frac{|f_\omega(x)-f_\omega(y)|}{D_{P,\alpha,\sigma}(x,y)} \leq
\frac{2|f_\omega(x)-f_\omega(y)|}{\|x-y\|},
\]
and the latter is maximized by finding two points $x$ and $y$ as close
together with nonzero numerator. In this case, $\|x-y\| = \varepsilon$
and $|f_\omega(x)-f_\omega(y)| = L \varepsilon/8$. Hence,
$|f_\omega(x)-f_\omega(y)|\leq L D_{P,\alpha,\sigma}(x,y)$ as required.
Now we show that each $P=P_\omega$ satisfies
\[
{\cal N}(P,\alpha,\sigma,\varepsilon) \leq \biggl(\frac{C}{\varepsilon
}
\biggr)^\xi
\]
for all $\varepsilon\geq n^{-{1}/({2+\xi})}$. Cover the top $A_1$ and
bottom $A_0$ of the cubes with Euclidean spheres of radius $\delta$.
There are $O((1/\delta)^{d-1})$ such spheres. The $D_{P,\alpha,\sigma}$
radius of each sphere is at most $\delta e^{-\alpha K(0)/\sigma^d}$.
Thus, these form an $\varepsilon$ covering as long as $\delta
e^{-\alpha K(0)/\sigma^d}\leq\varepsilon$. Thus the covering number of
the top and bottom is at most $2(1/\delta)^{d-1} \leq 2(1/(e^{\alpha
K(0)/\sigma^d}\varepsilon))^{d-1}$. Now cover the tendris with
one-dimensional segments of length $\delta$. The $D_{P,\alpha,\sigma}$
radius of each segment is at most $\delta e^{-\alpha/\sigma^d}$. Thus,
these form an $\varepsilon$ covering as long as $\delta e^{-\alpha
K(0)/\sigma^d}\leq\varepsilon$. Thus the covering number of the
tendrils is at most $q/\delta= n/\delta\leq n/(\varepsilon e^{\alpha
K(0)/\sigma^d})$. Thus we can cover the support with
\[
N(\varepsilon)\leq2 \biggl(\frac{1}{e^{\alpha K(0)/\sigma^d}\varepsilon
} \biggr)^{d-1} +
\frac{n}{\varepsilon e^{\alpha K(0)/\sigma^d}}
\]
balls of size $\varepsilon$. It follows from (\ref{eqalp}) that
$N(\varepsilon)\leq(1/\varepsilon)^\xi$ for $\varepsilon\geq
n^{-{1}/({2+\xi})}$ as required.\vspace*{9pt}

\textit{Claim} 2: For any $\omega$, and any $g \geq0$, $\int g(x)
\,dP_\omega(x) \geq \frac{1}{2} \int g(x) \,d\mu(x)$.\vspace*{9pt}

\textit{Proof}:
We have
\begin{eqnarray*}
\int_{S_\omega} g \,dP_\omega
&\geq& \int
_{C} g \,dP_\omega= \frac{1}{4}\int
_{A_0} g \,d\mu_0 + \frac{1}{4}\int
_{A_1} g \,d\mu_1 + \frac{\sum_j \int_{T_j} g \,d\nu_{j,\omega}}{2 q (1-\varepsilon)}
\\
&=& \frac{1}{4}\int_{A_0} g \,d\mu_0 +
\frac{1}{4}\int_{A_1} g \,d\mu_1 \\
&&{} +
\frac{(({1-2\varepsilon})\sum_j \int_{T_j} g \,d\nu_{j})/({1-\varepsilon})} {
2 q (1-2\varepsilon)} \times
\frac{{1}/{2} + q (1-2\varepsilon)}{ {1}/{2} + q (1-\varepsilon
)}
\\
&\geq& \frac{1}{2} \biggl( \frac{1}{4}\int_{A_0}
g \,d\mu_0 + \frac{1}{4}\int_{A_1} g d
\mu_1 + \frac{\sum_j \int_{T_j} g \,d\nu_{j}} {
2 q (1-2\varepsilon)} \biggr) = \frac{1}{2}\int g \,d\mu.
\end{eqnarray*}

\textit{Claim} 3:
For any $\omega,\nu\in\Omega$,
\[
d^2(f_\omega,f_\nu) = \frac{ \rho(\omega,\nu) L^2 \varepsilon^2 (1-2\varepsilon)} {
2 q (1-2\varepsilon)}.
\]

\textit{Proof}: This follows from direct calculation.\vspace*{9pt}

\textit{Claim} 4: If $\rho(\omega,\nu)=1$, then $\|P_\omega^n \wedge
P_\nu^n\| \geq1/(16e)$.\vadjust{\goodbreak}

\textit{Proof}: Suppose that $\rho(\omega,\nu)=1$. $P_\omega$ and
$P_\nu$ are the same everywhere except $T_{j,0}\cup T_{j,1}$, where $j$
is the index where $\omega$ and $\nu$ differ (assume $\omega_j=0$ and
$\nu_j=1$). Define $A= T_{j,0} \times\{0\}$ and
$B=T_{j,1}\times\{L\varepsilon\}$. Note that $A\cap B=\varnothing$. So,
\[
P_\omega(T_{j,0}\cup T_{j,1})= P_\omega(A)=P_\nu(T_{j,0}
\cup T_{j,1})=P_\nu(B)=\frac{1-\varepsilon
}{2q(1-\varepsilon)}
\]
and
\begin{eqnarray*}
\mathsf{TV}(P_\omega,P_\nu)&=&\bigl|P_\omega(A)-P_\nu(A)\bigr|
=\bigl|P_\omega(B)-P_\nu(B)\bigr|
\\
&=&\frac{1-\varepsilon}{2q(1-\varepsilon)}=\frac{1}{2q} = \frac{\varepsilon^{d-1}}{2}.
\end{eqnarray*}
Thus,
\[
\bigl\|P_\omega^n\land P_\nu^n\bigr\| \geq
\tfrac{1}{8} \bigl(1-\mathsf{TV}(P_\omega,P_\nu)
\bigr)^{2n} \geq\tfrac{1}{8} \bigl(1- \varepsilon^{d-1}/2
\bigr)^{2n}.
\]
Since $\varepsilon=n^{-{1}/({d-1})}$, this implies that
\[
\bigl\|P_\omega^n\land P_\nu^n\bigr\| \geq
\frac{1}{8} \biggl(1-\frac
{1}{2n} \biggr)^{2n} \geq
\frac{1}{16e}
\]
for all large $n$.\vspace*{9pt}

\textit{Completion of the proof.}
Recall that $\varepsilon= n^{-{1}/({d-1})}$.
Combining Assouad's lemma
(see Lemma 3 
in the supplement)
with the above claims,
we have
\begin{eqnarray*}
R_S &=& \inf_{\hat f\in{\cal S}_n}\sup_{P\in{\cal P}_{n,\xi}}
R_P(\hat f) \geq \inf_{\hat f\in{\cal S}_n}\sup
_{P\in{\cal P}_\Omega} R_P(\hat f) \geq \frac{1}{2}\inf
_{\hat f}\max_{\omega\in\Omega} \mathbb{E}_{\omega}
\bigl[d^2(f_{\omega},\hat f)\bigr]
\\
&\geq& \frac{q}{16} \times\frac{(L/8)^2\varepsilon^2(1-2\varepsilon
)}{2q(1-2\varepsilon)} \times\frac{1}{16e} = C
\frac{q\varepsilon^2(1-2\varepsilon)}{2q(1-2\varepsilon)}
\\
&\geq& C \varepsilon^2 = C n^{-{2}/({d-1})}.
\end{eqnarray*}
\upqed\end{pf}

\subsection{Semisupervised upper bound}

Now we state the upper bound for this class.

\begin{theorem}
Let $h = (n e^{2(2-\xi)})^{-{1}/({2+\xi})}$.
Then
%
\begin{equation}
\sup_{P\in{\cal P}_{n}}R(\hat f_{h,\alpha,\sigma})\leq \biggl(
\frac{C}{n} \biggr)^{{2}/({2+\xi})}.
\end{equation}
\end{theorem}

\begin{pf}
This follows from (C2), (C3) and
Corollary~\ref{corollaryupper}.
\end{pf}

\subsection{Comparison of lower and upper bound}

Combining the last two theorems we have:\vadjust{\goodbreak}

\begin{corollary}
Under the conditions of the previous theorem,
and assuming that $d > \xi+ 3$,
%
\begin{equation}
\frac{R_{SS}}{R_{S}} \preceq \biggl(\frac{1}{n} \biggr)
^{{2(d-3-\xi)}/({ (2+\xi)(d-1)})}\to0
\end{equation}
as $n\to\infty$.
\end{corollary}

This establishes the effectiveness of semi-supervised inference
in the minimax sense.

\section{Adaptive semisupervised inference}
\label{secadap}

We have established a bound on the risk
of the density-sensitive semisupervised kernel estimator.
The bound is achieved by using an estimate $\hat D_{\alpha,\sigma}$
of the
density-sensitive distance.
However, this requires knowing the density-sensitive parameter $\alpha
$, along with
other parameters.
It is critical to choose $\alpha$ (and $h$) appropriately, otherwise
we might
incur a large error if the semisupervised assumption does not hold, or
holds with a different density sensitivity value $\alpha$.
We consider two methods for choosing the parameters.

The following result shows that we can adapt to the correct degree of
semisupervisedness if cross-validation is used to select the
appropriate $\alpha,\sigma$ and $h$. This implies that the estimator
gracefully degrades to a supervised learner if the semisupervised
assumption (sensitivity of regression function to marginal density)
does not hold ($\alpha= 0$).

For any $f$, define the risk
$R(f) = \E[(f(X)-Y)^2]$ and the excess risk
$\cE(f) = R(f) - R(f^*) = \E[(f(X)-f^*(X))^2]$
where $f^*$ is the true regression function.
Let
${\cal H}$ be a finite set of bandwidths,
let ${\cal A}$ be a finite set of values for $\alpha$
and let $\Sigma$ be a finite set of values for $\sigma$.
Let $\theta= (h,\alpha,\sigma)$,
$\Theta= {\cal H}\times{\cal A}\times\Sigma$
and $J = |\Theta|$.

Divide the data into training data $T$
and validation data $V$.
For notational simplicity, let both sets have size $n$.
Let
${\cal F} = \{\hat f^T_{\theta}\}_{\theta\in\Theta}$
denote the semisupervised kernel estimators trained on
data $T$ using $\theta\in\Theta$.
For each $\hat f_{\theta}^T\in{\cal F}$ let
\[
\hat R^V \bigl(\hat f^T_{\theta}\bigr) =
\frac{1}{n}\sum^n_{i=1}\bigl(\hat
f^T_{\theta}(X_i)-Y_i
\bigr)^2,
\]
where the sum is over $V$.
Let $Y_i = f(X_i) + \varepsilon_i$ with $\varepsilon_i \stackrel
{\mathrm{i.i.d.}}{\sim} {\cal N}(0,\sigma^2)$.
Also, we assume that
$|f(x)|, |\hat f^T_{\theta}(x)| \leq M$, where $M>0$ is a
constant.\footnote{Note that the estimator can always be truncated if necessary.}

\begin{theorem}
\label{thmcrossval}
Let ${\cal F} = \{\hat f^T_{\theta}\}_{\theta\in\Theta}$
denote the semisupervised kernel estimators trained on
data $T$ using $\theta\in\Theta$.
Use validation data $V$ to pick
\[
\hat\theta= \mathop{\arg\min}_{\theta\in\Theta} \hat R^V\bigl(\hat
f^T_{\theta}\bigr)\vadjust{\goodbreak}
\]
and define the corresponding estimator
$\hat f = \hat f_{\hat\theta}$. Then,
for every $0 < \delta< 1$,
%
\begin{equation}
\E\bigl[\cE(\hat f_{\theta})\bigr] \leq\frac{1}{1-a} \biggl[\min
_{\theta\in\Theta} \E\bigl[\cE(\hat f_{\theta})\bigr] +
\frac{\log(J)/\delta)}{nt} \biggr] + 4\delta M^2,
\end{equation}
where $0<a<1$ and $0<t < 15/(38(M^2+\sigma^2))$ are constants. $\E$
denotes expectation over
everything that is random.
\end{theorem}

\begin{pf}
First,\vspace*{1pt} we derive a general concentration of $\hat\cE(f)$ around $\cE(f)$
where
$\hat\cE(f) = \hat R(f) - \hat R(f^*) = -\frac1{n}\sum^n_{i=1}U_i$
and $U_i = -(Y_i-f(X_i))^2+(Y_i-f^*(X_i))^2$.

If the variables $U_i$ satisfy the following moment condition:
\[
\E\bigl[\bigl|U_i - \E[U_i]\bigr|^k\bigr] \leq
\frac{\mathsf{Var}(U_i)}{2}k! r^{k-2}
\]
for some $r>0$, then the Craig--Bernstein (CB) inequality [Craig (\citeyear{Cra33})]
states that with probability $>1-\delta$,
\[
\frac1{n}\sum^n_{i=1}
\bigl(U_i - \E[U_i]\bigr) \leq\frac{\log(1/\delta
)}{nt} +
\frac{t \operatorname{\mathsf{Var}}(U_i)}{2(1-c)}
\]
for $0\leq tr \leq c <1$. The moment conditions are satisfied by
bounded random variables as well as Gaussian random variables; see, for example,
\citet{RandprojJHaupt}.

To apply this inequality, we first show that $\mathsf{Var}(U_i) \leq
4(M^2+\sigma^2)\cE(f)$ since
$Y_i = f(X_i) + \varepsilon_i$ with $\varepsilon_i \stackrel{\mathrm{i.i.d.}}{\sim}
{\cal N}(0,\sigma^2)$.
Also, we assume that $|f(x)|$, $|\hat f(x)| \leq M$, where $M>0$ is a constant.
\begin{eqnarray*}
\mathsf{Var}(U_i) &\leq& \E\bigl[U_i^2\bigr]
= \E\bigl[\bigl(-\bigl(Y_i-f(X_i)\bigr)^2+
\bigl(Y_i-f^*(X_i)\bigr)^2
\bigr)^2\bigr]
\\
&=& \E\bigl[\bigl(-\bigl(f^*(X_i)+\varepsilon_i-f(X_i)
\bigr)^2+(\varepsilon_i)^2\bigr)^2
\bigr]
\\
&=& \E\bigl[\bigl(-\bigl(f^*(X_i)-f(X_i)
\bigr)^2-2\varepsilon_i\bigl(f^*(X_i) -
f(X_i)\bigr)\bigr)^2\bigr]
\\
&\leq& 4 M^2 \cE(f) + 4 \sigma^2 \cE(f) = 4
\bigl(M^2+\sigma^2\bigr)\cE(f).
\end{eqnarray*}

Therefore using the CB inequality we get, with probability
$>1-\delta$,
\[
\cE(f) - \hat\cE(f) \leq\frac{\log(1/\delta)}{nt} + \frac{t
2(M^2+\sigma^2)\cE(f)}{(1-c)}.
\]
Now set $c = tr = 8t(M^2+\sigma^2)/15$ and let $t < 15/(38(M^2+\sigma
^2))$. With this choice,
$c < 1$ and define
\[
a = \frac{t 2(M^2+\sigma^2)}{(1-c)} <1.
\]
Then, using $a$ and rearranging terms, with probability $>1-\delta$,
\[
(1-a)\cE(f) - \hat\cE(f) \leq\frac{\log(1/\delta)}{nt},
\]
where $t < 15/(38(M^2+\sigma^2))$.

Then, using the previous
concentration result, and taking union bound over all $f\in{\cal F}$,
we have with probability $>1-\delta$,
\[
\cE(f) \leq\frac1{1-a} \biggl[\hat\cE^V(f) + \frac{\log(J/\delta
)}{nt}
\biggr].
\]
Now,
\begin{eqnarray*}
\cE(\hat f_{\hat\theta})
&=& R(\hat f_{\hat\theta
}) - R
\bigl(f^*\bigr)
\\
&\leq& \frac1{1-a} \biggl[\hat R^V(\hat
f_{\hat\theta}) - \hat R^V\bigl(f^*\bigr) + \frac{\log(J/\delta)}{nt}
\biggr]
\\
&\leq& \frac1{1-a} \biggl[\hat R^V(f) - \hat
R^V\bigl(f^*\bigr) + \frac{\log(J/\delta)}{nt} \biggr].
\end{eqnarray*}
Taking expectation with respect to validation dataset,
\[
\E_V\bigl[\cE(\hat f_{\hat\theta})\bigr] \leq\frac1{1-a}
\biggl[R(f) - R\bigl(f^*\bigr) + \frac{\log(J/\delta)}{nt} \biggr] + 4\delta
M^2.
\]
Now taking expectation with respect to training dataset,
\[
\E_{\mathrm{TV}}\bigl[\cE(\hat f_{\hat\theta})\bigr] \leq \frac1{1-a}
\biggl[\E_T\bigl[R(f) - R\bigl(f^*\bigr)\bigr] + \frac{\log(J/\delta
)}{nt}
\biggr] + 4 \delta M^2.
\]
Since this holds for all $f \in{\cal F}$, we get
\[
\E_{\mathrm{TV}}\bigl[\cE(\hat f_{\hat\theta})\bigr] \leq \frac1{1-a}
\biggl[\min_{f\in{\cal F}}\E_T\bigl[\cE(f)\bigr] +
\frac{\log(J/\delta)}{nt} \biggr] + 4 \delta M^2.
\]
The result follows.
\end{pf}

In practice, both
$\Theta$ may be taken to be
of size $n^a$ for some $a>0$.
Then we can approximate the optimal $h,\sigma$ and $\alpha$
with sufficient accuracy to achieve the optimal rate.
Setting $\delta= 1/(4 M^2n)$, we then see that
the penalty for
adaptation is
$\frac{\log(J/\delta)}{nt} + \delta M = O(\log n /n)$
and hence introduces only a logarithmic term.

\begin{Remark*}
Cross-validation is not the only way to adapt. For example, the
adaptive method in \citet{Kpotufe} can also be used here.
\end{Remark*}

\section{Simulation results}
\label{secsim}

\begin{figure}

\includegraphics{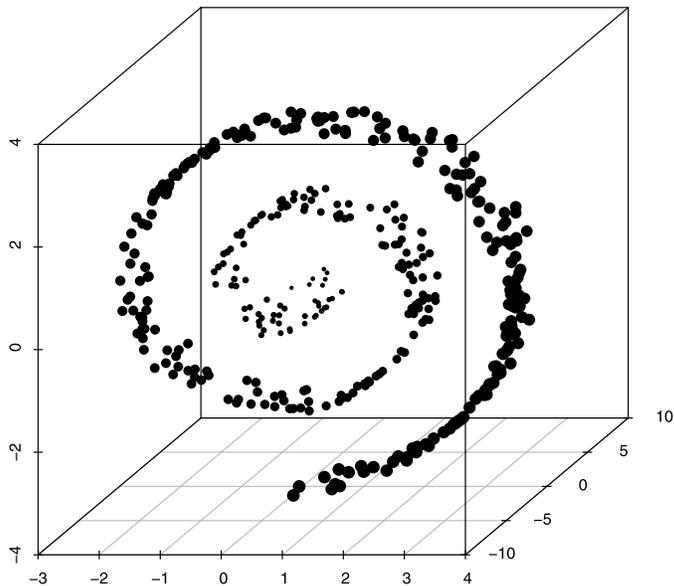}

\caption{The swiss roll data set. Point size represents regression function.}
\label{figswissroll}
\end{figure}

In this section we describe the results of a series of numerical
experiments on a simulated data\vadjust{\goodbreak} set to demonstrate the effect of using
the exponential version of the density sensitive metric for small, labeled
sample sizes. For the marginal distribution of $X$, we used a
slightly modified version of the swiss roll distribution used in
\citet{Culp2011SSS}. Figure~\ref{figswissroll} shows a sample from
this distribution, where the point size represents the response $Y$.
We repeatedly sampled $N=400$ points from this
distribution, and computed the mean squared error of the kernel
regression estimator using a set of values for $\alpha$ and for
labeled sample size ranging from $n=5$ to $n=320$.
We used the approximation method described in the supplement [see
equation (10)] 
with the number of nearest neighbors used set to $k=20$.

\begin{figure}

\includegraphics{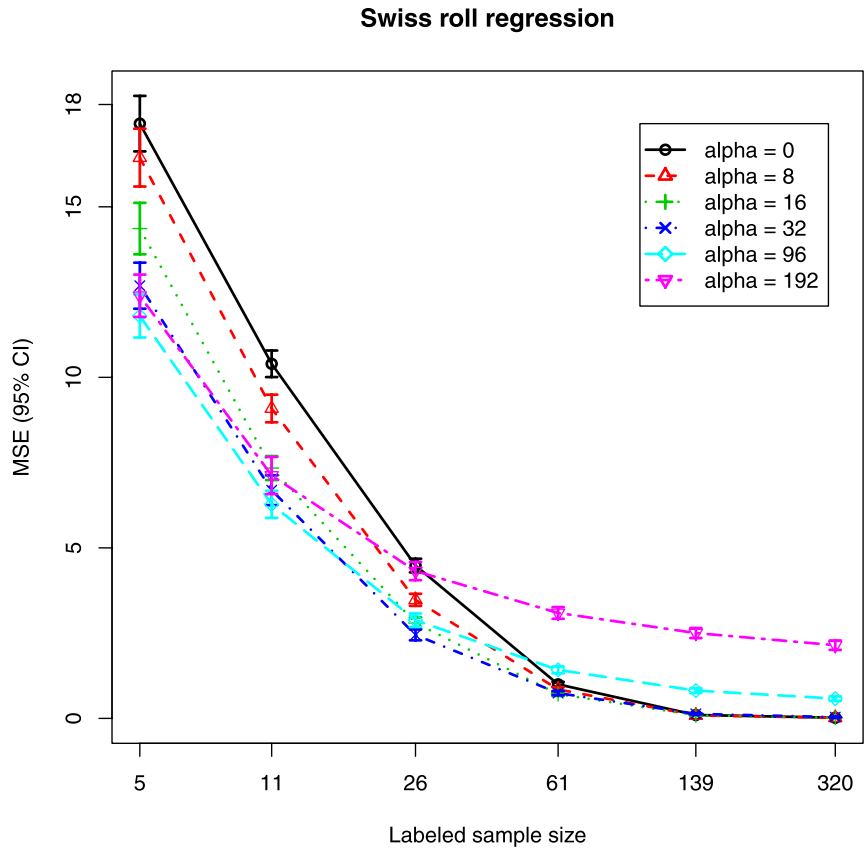}

\caption{MSE of kernel regression on the swiss roll data set for a
range of labeled sample sizes using different values of $\alpha$.}
\label{figsimulation}
\end{figure}

Figure~\ref{figsimulation} shows the average results after $300$
repetitions of this procedure with error bars indicating a $95\%$
confidence interval. As expected, we observe that for small labeled
sample sizes, increasing $\alpha$ can decrease the error. But as the
labeled sample size increases, using the density sensitive metric
becomes decreasingly beneficial, and can even hurt.

\section{Discussion}
\label{secdisc}

Semisupervised methods are very powerful, but like all methods, they
only work under certain conditions. We have shown that, under certain
conditions, semisupervised methods provably outperform supervised
methods. In particular, the advantage of semisupervised methods is
mainly when the distribution $P_X$ of $X$ is concentrated near a
low-dimensional set rather than when $P_X$ is smooth.

We introduced a family of estimators indexed by a parameter $\alpha$.
This parameter controls the strength of the semi-supervised assumption.
The behavior of the semi-supervised method depends critically on
$\alpha$. Finally, we showed that cross-validation can be used to
automatically adapt to $\alpha$ so that $\alpha$ does not need to be
known. Hence, our method takes advantage of the unlabeled data when the
semi-supervised assumption holds, but does not add extra bias when the
assumption fails. Our simulations confirm that our proposed estimator
has good risk when the semi-supervised smoothness holds.\looseness=-1

The analysis in this paper can be extended in several ways. First, it
is possible to use other density sensitive metrics such as the
diffusion distance [\citet{wasserman08spectral}]. Second, we
defined a method to estimate the density sensitive metric that works
under broader conditions than the two existing methods due to
\citet{orlitsky} and \citet{bigral11}. We suspect that faster
methods can be developed. Finally, other estimators besides kernel
estimators can be used. We will report on these extensions elsewhere.


\begin{supplement}
\stitle{Supplement to ``Density-sensitive semisupervised inference''\\}
\slink[doi]{10.1214/13-AOS1092SUPP} 
\sdatatype{.pdf}
\sfilename{aos1092\_supp.pdf}
\sdescription{Contains technical details, proofs and extensions.}
\end{supplement}


\printaddresses

\end{document}